# On the Tightness of Convex Optimal Power Flow Model Based on Power Loss Relaxation

Zhao Yuan
Electrical Power Systems Laboratory (EPS-Lab)
University of Iceland, Reykjavik, Iceland
zhaoyuan@hi.is

*Abstract*—Optimal power flow (OPF) is the fundamental mathematical model in power system operations. Improving the solution quality of OPF provide huge economic and engineering benefits. The convex reformulation of the original nonconvex alternating current OPF (ACOPF) model gives an efficient way to find the global optimal solution of ACOPF but suffers from the relaxation gaps. The existence of relaxation gaps hinders the practical application of convex OPF due to the AC-infeasibility problem. We evaluate and improve the tightness of the convex ACOPF model in this paper. Various power networks and nodal loads are considered in the evaluation. A unified evaluation framework is implemented in Julia programming language. This evaluation shows the sensitivity of the relaxation gap and helps to benchmark the proposed tightness reinforcement approach (TRA). The proposed TRA is based on the penalty function method which penalizes the power loss relaxation in the objective function of the convex ACOPF model. A heuristic penalty algorithm is proposed to find the proper penalty parameter of the TRA. Numerical results show relaxation gaps exist in test cases especially for large-scale power networks under low nodal power loads. TRA is effective to reduce the relaxation gap of the convex ACOPF model.

## Nomenclature

**Sets:**
$\mathcal{N}$    Nodes (or buses).
$\mathcal{L}$    Lines (or branches).

**Variables:**
$f(p_n)$    Original objective function.
$f^M(p_n)$    Modified objective function.
$p_n, q_n$    Active and reactive power generation at node $n$.
$p_{s_l}, q_{s_l}$    Sending-end active and reactive power flow of line $l$.
$p_{r_l}, q_{r_l}$    Receiving-end active and reactive power flow of line $l$.
$p_{o_l}, q_{o_l}$    Active and reactive power loss of line $l$.
$v_n, V_n$    Phase-to-ground voltage magnitude and voltage square at node $n$.
$v_{s_l}, v_{r_l}$    Sending- and receiving- end phase-to-ground voltage of line $l$.
$V_{s_l}, V_{r_l}$    Sending- and receiving- end phase-to-ground voltage square of line $l$.
$\theta_n$    Phase-to-ground voltage phase angle at node $n$.
$\theta_l$    Phase angle difference between the sending- and receiving- end phase-to-ground voltages of line $l$.
$\theta_{s_l}, \theta_{r_l}$    Phase angles of sending- and receiving- end phase-to-ground voltage of line $l$.
$K_{o_l}$    Equivalent ampacity constraint for power loss of line $l$.

**Parameters:**
$A^+_{nl}, A^-_{nl}$    Node to line (branch) incidence matrix.
$X_l, R_l$    Longitudinal reactance and resistance of line $l$ modelled as a passive Π-model.
$G_n, B_n$    Shunt conductance and susceptance of node $n$.
$B_{s_l}, B_{r_l}$    Sending- and receiving -end shunt susceptance of line $l$.
$\widetilde{K}_l, K_l$    Actual and approximated ampacity of line $l$.
$p_n^{min}, p_n^{max}$    Lower and upper bound of $p_n$.
$q_n^{min}, q_n^{min}$    Lower and upper bound of $q_n$.
$\theta_l^{min}, \theta_l^{min}$    Lower and upper bound of $\theta_l$.
$V_n^{min}, V_n^{min}$    Lower and upper bound of $V_n$.
$p_{d_n}, q_{d_n}$    Active and reactive power load of node $n$.
$\alpha_n, \beta_n, \gamma_n$    Cost parameters of active power generation.
$\xi, \Delta\xi$    Penalty coefficient and increment.

## I. Introduction

Efficient power system operations rely on accurate and adequate mathematical models of power networks [1]. Optimal power flow (OPF), as one of the most widely used mathematical optimization models, is executed frequently from every several minutes to hours in the energy management system (EMS) or supervisory control and data acquisition (SCADA) of power control centers around the world [2]–[4]. To name a few, the applications of OPF range from economic dispatch, transmission network expansion, security assessment, reliability management to multi-energy-carriers network operations [5]–[10]. The original OPF model considers the full alternating current (AC) circuit laws of three-phase power system is called ACOPF. Huge research efforts have put in the research of ACOPF in terms of its solvebility, approximation, convexity and computational complexity. The nonconvex nature of ACOPF makes it hard to find the global optimal solution. Recent advances in convex OPF based on power loss relaxation are promising approaches to improve the solution efficiency towards global optimality [11]–[14]. Authors in [14] rigorously prove the existence and uniqueness of global optimal solution of the second-order cone (SOC) based ACOPF model (SOC-ACOPF). It is also shown in [14] that relaxation gaps exist in some test cases. In order to apply the solution of SOC-ACOPF in actual power system operations, it is necessary to ensure the AC-feasibility of the solution i.e. to make sure that the solutions from the SOC-ACOPF model satisfy the original ACOPF constraints. Motivated in this aspect, an extensive investigation of the relaxation gap of the SOC-ACOPF model is firstly conducted in this paper. We then propose to use penalty function to improve the AC-feasibility of the solution from the SOC-ACOPF model. Though approaches including bound tightening approach and the sequential programming method are proposed and investigated in [6], [15] to find feasible solutions of the ACOPF models, no previous research work have applied penalty function approaches to the SOC-ACOPF model in [14] which is the focus of this work.

In this paper, the relaxation gap of SOC-ACOPF model is evaluated in a unified framework implemented in Julia programming language [16]. We also consider various power networks and power load levels in the evaluations. In order to guarantee the AC-feasibility of the SOC-ACOPF model, we propose to use penalty function method in the TRA. The proper penalty function parameter can be found by using the proposed heuristic penalty algorithm. Extensive numerical validations and discussions are provided. The rest of this paper is organized as follows. Section II presents the SOC-ACOPF model. Section III explains the relaxation gap metrics as the quantitative evaluation measures. Section IV proposes the TRAs. Section V gives the numerical results of the relaxation tightness evaluation and the performance of TRAs. Section VI concludes.

## II. Convex OPF Model Based on Power Loss Relaxation

The convex OPF model is presented in (1). This formulation is based on our recent work in [14]. Since second-order cone (SOC) is the key mathematical method to derive this model, we denote this model as SOC-ACOPF in this paper. The convexity and accuracy of this model have been rigorously proved in [14].

$$\text{Minimize} \quad f(p_n) = \sum_n (\alpha_n p_n^2 + \beta_n p_n + \gamma_n) \tag{1a}$$

subject to

$$p_n - p_{d_n} = \sum_l (A_{nl}^+ p_{s_l} - A_{nl}^- p_{o_l}) + G_n V_n, \ \forall n \in \mathcal{N} \tag{1b}$$

$$q_n - q_{d_n} = \sum_l (A_{nl}^+ q_{s_l} - A_{nl}^- q_{o_l}) - B_n V_n, \ \forall n \in \mathcal{N} \tag{1c}$$

$$V_{s_l} - V_{r_l} = 2R_l p_{s_l} + 2X_l q_{s_l} - R_l p_{o_l} - X_l q_{o_l}, \ \forall l \in \mathcal{L} \tag{1d}$$

$$\theta_l = X_l p_{s_l} - R_l q_{s_l}, \ \forall l \in \mathcal{L} \tag{1e}$$

$$K_{o_l} \geq q_{o_l} \geq \frac{p_{s(r)_l}^2 + q_{s(r)_l}^2}{V_{s(r)_l}} X_l, \ \forall l \in \mathcal{L} \tag{1f}$$

$$p_{o_l} X_l = q_{o_l} R_l, \ \forall l \in \mathcal{L} \tag{1g}$$

$$\theta_l = \theta_{s_l} - \theta_{r_l}, \ \forall l \in \mathcal{L} \tag{1h}$$

$$V_{s_l} V_{r_l} \sin^2(\theta_l^{max}) \geq \theta_l^2, \ \forall l \in \mathcal{L} \tag{1i}$$

$$V_n \in (V_n^{min}, V_n^{max}), \ \forall n \in \mathcal{N} \tag{1j}$$

$$\theta_l \in (\theta_l^{min}, \theta_l^{max}), \ \forall l \in \mathcal{L} \tag{1k}$$

$$\theta_n \in (\theta_n^{min}, \theta_n^{max}), \ \forall n \in \mathcal{N} \tag{1l}$$

$$p_n \in (p_n^{min}, p_n^{max}), \ \forall n \in \mathcal{N} \tag{1m}$$

$$q_n \in (q_n^{min}, q_n^{max}), \ \forall n \in \mathcal{N} \tag{1n}$$

Where $n \in \mathcal{N}$ is the index of nodes. $l \in \mathcal{L}$ is the index of lines. We use the expression of economic cost of power generation as the objective function $f(p_n)$ in this paper. Since the objective function is quadratic, it is convex. $(\alpha_n, \beta_n, \gamma_n) \geq 0$ are the cost parameters of the nodal active power generation. Equations (1b) and (1c) are to represent the active and reactive power balance constraints. $G_n, B_n$ are the shunt conductance and susceptance of node $n$. $A_{nl}^+$ and $A_{nl}^-$ are the node-to-branch incidence matrices of the power network. Specifically, $A_{nl}^+ = 1$, $A_{nl}^- = 0$ if $n$ is the sending-end of branch $l$, and $A_{nl}^+ = -1$, $A_{nl}^- = -1$ if $n$ is the receiving-end of branch $l$. Variables $p_n, q_n$ are the active and reactive power generations at node $n$. Parameters $p_{d_n}, q_{d_n}$ are the active and reactive power loads of node $n$. Variable $V_n = v_n^2$ is the phase-to-ground voltage magnitude square at node $n$. Note the reason of using $V_n$ (instead of $v_n$) is to guarantee the convexity of the SOC-ACOPF model. It is quite straightforward to recover the value of $v_n$ by $v_n = \sqrt{V_n}$ after solving the SOC-ACOPF model. Variables $p_{o_l}, q_{o_l}$ are the active and reactive power loss of branch $l$. Equations (1d)-(1e) are derived by taking the magnitude and phase angle of the voltage drop phasor along branch $l$ respectively. $v_{s_l}, v_{r_l}$ are the sending- and receiving- end phase-to-ground voltage of branch $l$. $V_{s(r)_l} = v_{s(r)_l}^2$ are voltage magnitude squares. $\theta_l = \theta_{s_l} - \theta_{r_l}$ is the phase angle difference between the sending- and receiving- end voltages of branch $l$. $\theta_{s_l}, \theta_{r_l}$ are the phase angles of sending- and receiving- end phase-to-ground voltages of branch $l$. To guarantee this derivation is valid, we assume $(\theta_l^{min}, \theta_l^{max}) \subseteq (-\frac{\pi}{2}, \frac{\pi}{2})$. This assumption is valid in power system operations. Constraint (1f) is the reactive power loss relaxation. The tightness of this relaxation is the focus of this paper. When the relaxation is tight, the solution of the SOC-ACOPF model is AC-feasible and thus global optimal for the original ACOPF model. Otherwise, the solution is not feasible for the original ACOPF model. In this case, a feasible solution recovery procedure is required to recover a AC-feasible solution. Equation (1g) represents the active and reactive power loss relationship. $\widetilde{K}_{ol}$ is the upper bound of reactive power loss which is used to constrain the capacity of branch $l$ equivalently. Constraint (1i) is used to improve the AC feasibility of the SOC-ACOPF model. Constraints (1j)-(1n) are bounds for voltage magnitude, voltage phase angle difference, nodal active power injection and nodal reactive power generation variables. $p_n^{min}, p_n^{max}$ are the lower and upper bounds of $p_n$. $q_n^{min}, q_n^{max}$ are the lower and upper bounds of $q_n$. $\theta_l^{min}, \theta_l^{max}$ are the lower and upper bounds of $\theta_l$. This model is valid for both radial and meshed power networks. The upper bound of reactive power loss $K_{ol}$ is quantified as (1o) [14]:

$$K_{o_l} = (\widetilde{K}_l - V_{s(r)_l} B_{s(r)_l}^2 + 2q_{s(r)_l} B_{s(r)_l}) X_l \tag{1o}$$

**TABLE I.** Relaxation Gaps of Active Power Loss

| Nodal Load | case9 | IEEE14 | case30 | IEEE57 | case89pegase | IEEE118 | ACTIVSg200 | IEEE300 | ACTIVSg500 |
|---|---|---|---|---|---|---|---|---|---|
| 5% | 3.24E-02 | 0.00E+00 | 0.00E+00 | 0.00E+00 | 2.06E+00 | 0.00E+00 | 3.32E+00 | 1.46E-02 | 1.69E+01 |
| 10% | 0.00E+00 | 0.00E+00 | 0.00E+00 | 0.00E+00 | 1.23E+00 | 0.00E+00 | 3.23E+00 | 1.15E-02 | 1.31E+01 |
| 15% | 0.00E+00 | 0.00E+00 | 0.00E+00 | 0.00E+00 | 1.04E+00 | 0.00E+00 | 3.04E+00 | 6.62E-03 | 9.35E+00 |
| 20% | 0.00E+00 | 0.00E+00 | 0.00E+00 | 0.00E+00 | 5.99E-02 | 0.00E+00 | 2.73E+00 | 9.95E-04 | 5.57E+00 |
| 25% | 0.00E+00 | 0.00E+00 | 0.00E+00 | 0.00E+00 | 7.98E-03 | 0.00E+00 | 2.46E+00 | 3.79E-13 | 1.80E+00 |
| 30% | 0.00E+00 | 0.00E+00 | 0.00E+00 | 0.00E+00 | 2.87E-11 | 0.00E+00 | 2.19E+00 | 3.48E-13 | 4.19E-02 |
| 35% | 0.00E+00 | 0.00E+00 | 0.00E+00 | 0.00E+00 | 2.83E-11 | 0.00E+00 | 2.12E+00 | 3.16E-13 | 2.20E-02 |
| 40% | 0.00E+00 | 0.00E+00 | 0.00E+00 | 0.00E+00 | 2.94E-11 | 0.00E+00 | 1.84E+00 | 2.87E-13 | 2.08E-02 |
| 45% | 0.00E+00 | 0.00E+00 | 0.00E+00 | 0.00E+00 | 3.14E-11 | 0.00E+00 | 1.61E+00 | 2.60E-13 | 2.89E-02 |
| 50% | 0.00E+00 | 0.00E+00 | 0.00E+00 | 0.00E+00 | 3.34E-11 | 0.00E+00 | 1.45E+00 | 2.34E-13 | 1.93E-02 |
| 55% | 0.00E+00 | 0.00E+00 | 0.00E+00 | 0.00E+00 | 3.57E-11 | 0.00E+00 | 1.35E+00 | 0.00E+00 | 1.84E-02 |
| 60% | 0.00E+00 | 0.00E+00 | 0.00E+00 | 0.00E+00 | 3.85E-11 | 0.00E+00 | 9.10E-01 | 1.85E-13 | 1.69E-02 |
| 65% | 0.00E+00 | 0.00E+00 | 0.00E+00 | 0.00E+00 | 4.53E-11 | 0.00E+00 | 6.35E-01 | 1.62E-13 | 4.47E-03 |
| 70% | 0.00E+00 | 0.00E+00 | 0.00E+00 | 0.00E+00 | 5.55E-11 | 0.00E+00 | 3.34E-01 | 1.40E-13 | 6.19E-12 |
| 75% | 0.00E+00 | 0.00E+00 | 0.00E+00 | 0.00E+00 | 7.36E-11 | 0.00E+00 | 4.05E-02 | 1.15E-13 | 2.46E-11 |
| 80% | 0.00E+00 | 0.00E+00 | 0.00E+00 | 0.00E+00 | 1.13E-10 | 0.00E+00 | 6.35E-03 | 9.32E-14 | 1.49E-06 |
| 85% | 0.00E+00 | 0.00E+00 | 0.00E+00 | 0.00E+00 | 2.55E-10 | 0.00E+00 | 6.26E-04 | 6.99E-14 | 7.38E-04 |
| 90% | 0.00E+00 | 0.00E+00 | 0.00E+00 | 0.00E+00 | 2.48E-03 | 0.00E+00 | 3.86E-12 | 5.64E-14 | 4.22E-04 |
| 95% | 0.00E+00 | 0.00E+00 | 0.00E+00 | 0.00E+00 | 1.18E-02 | 0.00E+00 | 3.64E-12 | 4.57E-14 | 1.75E-13 |
| 100% | 0.00E+00 | 0.00E+00 | 0.00E+00 | 0.00E+00 | 1.85E-02 | 0.00E+00 | -3.78E-13 | 3.39E-14 | 2.95E-15 |

**TABLE II.** Relaxation Gaps of Reactive Power Loss

| Nodal Load | case9 | IEEE14 | case30 | IEEE57 | case89pegase | IEEE118 | ACTIVSg200 | IEEE300 | ACTIVSg500 |
|---|---|---|---|---|---|---|---|---|---|
| 5% | 2.62E-01 | 3.14E-01 | 2.00E-02 | 2.42E-01 | 5.67E+00 | 2.71E+00 | 7.71E-01 | 7.58E+00 | 1.30E+00 |
| 10% | 3.82E-01 | 2.47E-01 | 3.05E-02 | 2.25E-01 | 3.37E+00 | 2.98E+00 | 7.71E-01 | 6.25E+00 | 1.30E+00 |
| 15% | 3.82E-01 | 1.77E-01 | 3.05E-02 | 2.13E-01 | 1.77E+00 | 3.16E+00 | 7.71E-01 | 5.74E+00 | 1.30E+00 |
| 20% | 2.39E-01 | 8.96E-02 | 3.03E-02 | 1.92E-01 | 6.74E-01 | 3.08E+00 | 7.71E-01 | 5.23E+00 | 1.30E+00 |
| 25% | 2.35E-01 | 1.56E-01 | 2.99E-02 | 1.69E-01 | 1.06E+00 | 2.97E+00 | 7.70E-01 | 4.55E+00 | 1.30E+00 |
| 30% | 2.31E-01 | 1.43E-01 | 2.96E-02 | 1.46E-01 | 7.84E-01 | 2.86E+00 | 7.70E-01 | 4.21E+00 | 4.02E-01 |
| 35% | 2.27E-01 | 1.30E-01 | 2.92E-02 | 1.21E-01 | 6.63E-01 | 2.90E+00 | 7.70E-01 | 3.92E+00 | 3.01E-01 |
| 40% | 2.22E-01 | 1.19E-01 | 2.87E-02 | 1.02E-01 | 5.82E-01 | 2.90E+00 | 7.69E-01 | 3.69E+00 | 2.95E-01 |
| 45% | 2.17E-01 | 1.07E-01 | 2.82E-02 | 1.17E-01 | 5.11E-01 | 3.10E+00 | 7.69E-01 | 3.39E+00 | 2.78E-01 |
| 50% | 2.13E-01 | 9.33E-02 | 2.78E-02 | 1.28E-01 | 4.42E-01 | 3.31E+00 | 7.69E-01 | 3.03E+00 | 2.83E-01 |
| 55% | 2.10E-01 | 6.97E-02 | 2.74E-02 | 1.37E-01 | 3.78E-01 | 3.49E+00 | 7.69E-01 | 2.65E+00 | 2.77E-01 |
| 60% | 2.14E-01 | 2.49E-02 | 2.69E-02 | 1.46E-01 | 3.17E-01 | 3.64E+00 | 7.68E-01 | 2.32E+00 | 2.70E-01 |
| 65% | 2.01E-01 | 1.39E-10 | 2.65E-02 | 1.57E-01 | 2.62E-01 | 3.80E+00 | 7.68E-01 | 2.05E+00 | 2.77E-01 |
| 70% | 2.03E-01 | 4.03E-10 | 2.60E-02 | 1.70E-01 | 2.28E-01 | 3.95E+00 | 7.68E-01 | 1.81E+00 | 9.09E-10 |
| 75% | 1.97E-01 | 4.04E-10 | 2.54E-02 | 1.72E-01 | 1.93E-01 | 3.99E+00 | 2.89E-01 | 1.60E+00 | 8.38E-10 |
| 80% | 1.90E-01 | 4.06E-10 | 2.47E-02 | 1.72E-01 | 1.59E-01 | 3.99E+00 | 1.95E-01 | 1.39E+00 | 5.05E-05 |
| 85% | 1.97E-01 | 4.08E-10 | 2.40E-02 | 1.53E-01 | 1.26E-01 | 3.99E+00 | 4.33E-02 | 1.20E+00 | 2.51E-02 |
| 90% | 1.91E-01 | 4.09E-10 | 2.32E-02 | 1.24E-01 | 1.97E-01 | 3.98E+00 | 2.82E-10 | 1.03E+00 | 1.44E-02 |
| 95% | 1.77E-01 | 4.10E-10 | 2.23E-02 | 9.71E-02 | 9.31E-01 | 3.98E+00 | 2.66E-10 | 8.30E-01 | 1.98E-11 |
| 100% | 1.71E-01 | 4.07E-10 | 2.00E-02 | 8.30E-02 | 1.47E+00 | 3.98E+00 | -1.90E-11 | 6.81E-01 | 3.35E-13 |

Where $\widetilde{K}_l$ is the actual ampacity of branch $l$ which is usually provided by the branch manufacturer, $B_{s(r)_l}$ is the shunt susceptance of branch $l$. For details of the derivations of (1o) according to the correct physical interpretation of the widely used transmission line $\Pi$-model, please refer to [14]. Note constraint (1o) is linear and thus convex.

## III. EVALUATION METRICS OF THE RELAXATION GAPS

For the SOC-ACOPF model-(1), we define the relaxation gap of active power loss $Gap_l^{po}$ as $Gap_l^{po}$ as $Gap_l^{po} := p_{o_l} - \frac{p_{s_l}^2 + q_{s_l}^2}{V_{s_l}} R_l, \forall l \in \mathcal{L}$ and the relaxation gap of reactive power loss $Gap_l^{qo}$ as $Gap_l^{qo} := q_{o_l} - \frac{p_{s_l}^2 + q_{s_l}^2}{V_{s_l}} X_l, \forall l \in \mathcal{L}$. The corresponding maximum relaxation gaps (of active and reactive power loss) are defined as $Gap^{po,max} := Maximum\{Gap_l^{po}, \forall l \in \mathcal{L}\}$ and $Gap^{qo,max} := Maximum\{Gap_l^{qo}, \forall l \in \mathcal{L}\}$. A fully AC-feasible solution of the SOC-ACOPF model-(1) means that $Gap^{po,max} = Gap^{qo,max} = 0$. When $Gap^{po,max} \neq 0$ or $Gap^{qo,max} \neq 0$, smaller values of $Gap^{po,max}, Gap^{qo,max}$ mean better solution quality in terms of AC-feasibility.

---

**Algorithm 1:** Heuristic Penalty Algorithm

**Result:** AC-feasible Solution of the SOC-ACOPF model-(1)

1 Initialization $\xi = \xi_0, k = 1$;
2 **do**
3     Solve the SOC-ACOPF model-(1);
4     $\xi = \xi + \Delta\xi$;
5     $k = k + 1$
6 **while** $Gap^{qo,max} > Gap^{qo,tol}$ (or $Gap^{po,max} > Gap^{po,tol}$) and $k < k^{max}$;

**TABLE III.** Tightened Relaxation Gaps of Active Power Loss

| Nodal Load | case9 | IEEE14 | case30 | IEEE57 | case89pegase | IEEE118 | ACTIVSg200 | IEEE300 | ACTIVSg500 |
|---|---|---|---|---|---|---|---|---|---|
| 5% | 8.75E-02 | 0.00E+00 | 0.00E+00 | 0.00E+00 | 2.21E+00 | 0.00E+00 | 3.33E+00 | 3.44E-02 | 1.69E+01 |
| 10% | 0.00E+00 | 0.00E+00 | 0.00E+00 | 0.00E+00 | 1.56E+00 | 0.00E+00 | 3.24E+00 | 3.38E-02 | 1.31E+01 |
| 15% | 0.00E+00 | 0.00E+00 | 0.00E+00 | 0.00E+00 | 3.69E-01 | 0.00E+00 | 3.03E+00 | 2.70E-02 | 9.35E+00 |
| 20% | 0.00E+00 | 0.00E+00 | 0.00E+00 | 0.00E+00 | 4.88E-01 | 0.00E+00 | 2.74E+00 | 1.55E-02 | 5.57E+00 |
| 25% | 0.00E+00 | 0.00E+00 | 0.00E+00 | 0.00E+00 | 9.63E-01 | 0.00E+00 | 2.47E+00 | 5.17E-03 | 1.80E+00 |
| 30% | 0.00E+00 | 0.00E+00 | 0.00E+00 | 0.00E+00 | 6.40E-12 | 0.00E+00 | 2.19E+00 | 5.84E-04 | 9.18E-02 |
| 35% | 0.00E+00 | 0.00E+00 | 0.00E+00 | 0.00E+00 | 6.36E-12 | 0.00E+00 | 2.12E+00 | 3.65E-04 | 3.00E-02 |
| 40% | 0.00E+00 | 0.00E+00 | 0.00E+00 | 0.00E+00 | 6.34E-12 | 0.00E+00 | 1.84E+00 | 9.78E-05 | 2.75E-02 |
| 45% | 0.00E+00 | 0.00E+00 | 0.00E+00 | 0.00E+00 | 1.98E-12 | 0.00E+00 | 1.62E+00 | 2.06E-13 | 1.58E-11 |
| 50% | 0.00E+00 | 0.00E+00 | 0.00E+00 | 0.00E+00 | 6.23E-12 | 0.00E+00 | 1.45E+00 | 1.80E-13 | 1.65E-11 |
| 55% | 0.00E+00 | 0.00E+00 | 0.00E+00 | 0.00E+00 | 6.28E-12 | 0.00E+00 | 1.35E+00 | 1.54E-13 | 1.71E-11 |
| 60% | 0.00E+00 | 0.00E+00 | 0.00E+00 | 0.00E+00 | 6.05E-12 | 0.00E+00 | 9.20E-01 | 1.31E-13 | 2.15E-11 |
| 65% | 0.00E+00 | 0.00E+00 | 0.00E+00 | 0.00E+00 | 6.18E-12 | 0.00E+00 | 6.32E-01 | 1.10E-13 | 2.22E-11 |
| 70% | 0.00E+00 | 0.00E+00 | 0.00E+00 | 0.00E+00 | 6.13E-12 | 0.00E+00 | 3.40E-01 | 0.00E+00 | 9.86E-13 |
| 75% | 0.00E+00 | 0.00E+00 | 0.00E+00 | 0.00E+00 | 6.16E-12 | 0.00E+00 | 9.92E-02 | 6.90E-14 | 6.63E-13 |
| 80% | 0.00E+00 | 0.00E+00 | 0.00E+00 | 0.00E+00 | 6.15E-12 | 0.00E+00 | 9.92E-02 | 4.93E-14 | 4.58E-13 |
| 85% | 0.00E+00 | 0.00E+00 | 0.00E+00 | 0.00E+00 | 6.13E-12 | 0.00E+00 | 3.54E-02 | 2.98E-14 | 4.33E-13 |
| 90% | 0.00E+00 | 0.00E+00 | 0.00E+00 | 0.00E+00 | 6.12E-12 | 0.00E+00 | 1.38E-12 | 2.13E-14 | 3.98E-04 |
| 95% | 0.00E+00 | 0.00E+00 | 0.00E+00 | 0.00E+00 | 5.87E-12 | 0.00E+00 | 7.83E-13 | 1.45E-14 | 4.42E-08 |
| 100% | 0.00E+00 | 0.00E+00 | 0.00E+00 | 0.00E+00 | 5.87E-12 | 0.00E+00 | 7.77E-13 | 2.48E-15 | 5.71E-13 |

**TABLE IV.** Tightened Relaxation Gaps of Active Power Loss

| Nodal Load | case9 | IEEE14 | case30 | IEEE57 | case89pegase | IEEE118 | ACTIVSg200 | IEEE300 | ACTIVSg500 |
|---|---|---|---|---|---|---|---|---|---|
| 5% | 3.82E-01 | 4.66E-10 | 1.64E-10 | 9.77E-03 | 5.44E+00 | 3.35E-11 | 7.71E-01 | 2.97E+00 | 1.30E+00 |
| 10% | 6.36E-10 | 4.12E-10 | 1.72E-10 | 6.28E-03 | 3.84E+00 | 3.44E-11 | 7.71E-01 | 2.67E+00 | 1.30E+00 |
| 15% | 6.25E-10 | 4.11E-10 | 1.76E-10 | 2.98E-03 | 9.09E-01 | 3.03E-11 | 7.71E-01 | 1.78E+00 | 1.30E+00 |
| 20% | 5.39E-10 | 4.10E-10 | 1.66E-10 | 3.46E-03 | 8.27E-01 | 3.03E-11 | 7.71E-01 | 7.32E-01 | 1.30E+00 |
| 25% | 4.23E-10 | 4.08E-10 | 1.67E-10 | 7.92E-03 | 1.63E+00 | 3.70E-11 | 7.70E-01 | 5.43E-01 | 1.30E+00 |
| 30% | 4.21E-10 | 4.07E-10 | 1.68E-10 | 8.78E-03 | 1.03E-10 | 3.04E-11 | 7.70E-01 | 2.48E-01 | 8.42E-01 |
| 35% | 4.20E-10 | 4.07E-10 | 1.69E-10 | 3.49E-03 | 1.02E-10 | 3.20E-11 | 7.70E-01 | 2.26E-01 | 2.75E-01 |
| 40% | 4.17E-10 | 4.05E-10 | 1.70E-10 | 1.41E-10 | 1.02E-10 | 1.79E-10 | 7.70E-01 | 2.03E-01 | 2.65E-01 |
| 45% | 4.14E-10 | 4.02E-10 | 1.70E-10 | 1.42E-10 | 3.58E-11 | 2.82E-10 | 7.69E-01 | 2.14E-01 | 7.28E-11 |
| 50% | 4.12E-10 | 4.01E-10 | 1.70E-10 | 1.42E-10 | 1.03E-10 | 1.81E-10 | 7.69E-01 | 2.31E-01 | 7.60E-11 |
| 55% | 4.11E-10 | 4.02E-10 | 1.71E-10 | 1.43E-10 | 1.04E-10 | 3.10E-11 | 7.69E-01 | 2.28E-01 | 7.91E-11 |
| 60% | 4.63E-10 | 3.99E-10 | 1.71E-10 | 1.43E-10 | 8.97E-11 | 3.12E-11 | 7.68E-01 | 2.13E-01 | 1.27E-10 |
| 65% | 4.61E-10 | 4.02E-10 | 1.71E-10 | 1.44E-10 | 8.86E-11 | 3.15E-11 | 7.68E-01 | 1.98E-01 | 1.28E-10 |
| 70% | 4.59E-10 | 4.04E-10 | 1.71E-10 | 1.44E-10 | 8.75E-11 | 3.17E-11 | 7.68E-01 | 1.85E-01 | 1.41E-11 |
| 75% | 4.57E-10 | 4.06E-10 | 1.71E-10 | 1.45E-10 | 8.65E-11 | 3.19E-11 | 6.40E-01 | 1.84E-01 | 9.50E-12 |
| 80% | 4.01E-10 | 4.08E-10 | 1.71E-10 | 1.45E-10 | 8.55E-11 | 3.22E-11 | 6.41E-01 | 1.83E-01 | 6.57E-12 |
| 85% | 4.53E-10 | 4.09E-10 | 1.71E-10 | 1.46E-10 | 8.44E-11 | 3.24E-11 | 2.29E-01 | 2.11E-01 | 6.21E-12 |
| 90% | 3.97E-10 | 4.10E-10 | 1.72E-10 | 1.46E-10 | 8.27E-11 | 3.25E-11 | 2.07E-11 | 2.88E-01 | 1.35E-02 |
| 95% | 4.48E-10 | 4.09E-10 | 1.73E-10 | 1.46E-10 | 8.04E-11 | 3.71E-11 | 9.79E-12 | 3.10E-01 | 2.08E-07 |
| 100% | 3.94E-10 | 4.09E-10 | 1.76E-10 | 1.60E-10 | 9.69E-11 | 3.25E-11 | 9.84E-12 | 2.66E-01 | 9.38E-12 |

## IV. TIGHTNESS REINFORCEMENT

In order to reinforce the tightness of relaxation in constraint (1f), we propose one heuristic approach based on penalty function. A penalty term $\xi \sum_l q_{ol}$ is added in the objective function $f(p_n)$ of the SOC-ACOPF model-(1) as: $f^M(p_n) = f(p_n) + \xi \sum_l q_{ol}, \forall l \in \mathcal{L}$. Where $f^M(p_n)$ is the modified objective function, $\xi > 0$ is the penalty coefficient parameter. Since we minimize the objective function, the added penalty term helps to reduce the relaxation gaps. This approach also gives the flexibility of modifying the penalty coefficient $\xi$ to reduce the relaxation gaps down to sufficient levels. It is worth to mention that since the original generation cost function is modified, the final solution of TRA is a compromise between optimality and AC-feasibility. In order to tune the penalty coefficient parameter $\xi$, we propose the heuristic penalty algorithm 1 to finally find the AC-feasible solution of the SOC-ACOPF model-(1). Where $0 < \xi_0 < 1$ is the initial penalty parameter, $k$ is the index of iteration, $k^{max}$ is the maximum allowed number of iterations, parameters $Gap^{po,tol}, Gap^{qo,tol}$ are the tolerances of active and reactive power loss relaxation gaps, parameter $0 < \Delta\xi << 1$ is used to increase the penalty parameter. This algorithm works iteratively to reduce the relaxation gaps and find a AC-feasible solution of the SOC-ACOPF model-(1). We show the performance of above two formulations of penalty functions in the numerical results Section.

## V. NUMERICAL RESULTS

The relaxation gap evaluation and TRAs are implemented in Julia programming language running on the 64-bit Windows 10 operating system. A personal computer with AMD Ryzen-7 2.9 GHz CPU and 16G RAM is deployed. The IPOPT solver in Julia is used to solve the SOC-ACOPF model and TRAs. We use the power network data from MATPOWER [17]. Evaluated power networks include case9, IEEE14, IEEE57, case89pegase, IEEE118, IEEE300 and ACTIVSg500 [18]–[21]. The banchmark nodal load or 100% nodal load is set same as the load level in the original data file from MATPOWER. We then modify the nodal power load from 5% to 100% in order to evaluate the relaxation gaps for different power load levels. The penalty coefficient is set as $\xi = 0.3$ in the numerical validations. The active power loss relaxation gap results are listed in Table I. For most cases, the relaxation gaps of active power loss are equal or very close to zero. There are some non-zero relaxation gaps of the active power loss when the nodal power loads are below 10% of the original power loads. The reactive power loss relaxation gap results are listed in Table II. These values are generally larger than the relaxation gaps of the active power loss. The reason is majorly due to the fact that reactance is much larger than resistance in transmission lines. We can also see that when the nodal power loads increase, the relaxation gaps decrease. However, the relaxation gaps of reactive power loss are still large which necessitates using TRA to tighten the relaxation gaps. The performance of penalty function based TRA is shown in Table III and Table IV. Compared with the relaxation gaps in Table I and Table II, the relaxation gaps are reduced. By using the test case of IEEE300, we show the performance of the heuristic penalty algorithm in Fig. 1. The relaxation gaps of both active power loss and reactive power loss are tightened to zero when the penalty coefficient parameter is larger or equal to one.

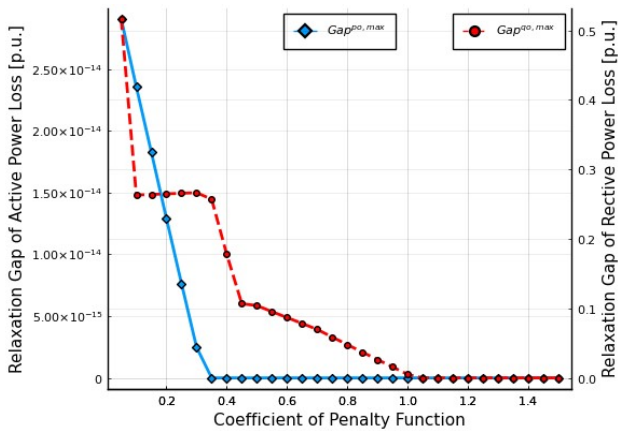

**Fig. 1.** Performance of the Heuristic Penalty Algorithm.

## VI. CONCLUSIONS

The contributions of this paper are two-fold: (1) extensively evaluating the relaxation gaps of the SOC-ACOPF model-(1) for various power networks and nodal loads; (2) proposing the TRA to improve the AC-feasibility of the solutions from the SOC-ACOPF model-(1). Numerical results show that the relaxation gap of the SOC-ACOPF model is sensitive to both power networks and nodal loads. The proposed TRA based on penalty functions is helpful to recover feasible solutions from the infeasible solutions of the SOC-ACOPF model-(1). Results in this paper can serve as benchmarks for future research of the SOC-ACOPF model-(1). Smart approaches for the penalty parameter selection of the TRA are expected for future works.

## REFERENCES

[1] A. Bergen and V. Vittal, *Power Systems Analysis*. Pearson/Prentice Hall, 2000.
[2] J. Carpentier, "Contribution to the economic dispatch problem," *Bulletin de la Societe Francoise des Electriciens*, vol. 3, no. 8, pp. 431–447, 1962.
[3] M. B. Cain, R. P. O?neill, A. Castillo *et al.*, "History of optimal power flow and formulations," *Federal Energy Regulatory Commission*, vol. 1, pp. 1–36, 2012.
[4] F. Capitanescu, J. M. Ramos, P. Panciatici, D. Kirschen, A. M. Marcolini, L. Platbrood, and L. Wehenkel, "State-of-the-art, challenges, and future trends in security constrained optimal power flow," *Electric Power Systems Research*, vol. 81, no. 8, pp. 1731–1741, 2011.
[5] W. Xu, A. Li, Y. Su, M. Zhu, X. Ouyang, and X. Yang, "Optimal expansion planning of ac/dc hybrid system integrated with vsc control strategy," in *2019 IEEE Innovative Smart Grid Technologies - Asia (ISGT Asia)*, 2019, pp. 3272–3276.
[6] Z. Yuan, M. R. Hesamzadeh, and D. Biggar, "Distribution locational marginal pricing by convexified ACOPF and hierarchical dispatch," *IEEE Transactions on Smart Grid*, vol. PP, no. 99, pp. 1–1, 2016.
[7] J. Condren and T. W. Gedra, "Expected-security-cost optimal power flow with small-signal stability constraints," *IEEE Transactions on Power Systems*, vol. 21, no. 4, pp. 1736–1743, 2006.
[8] Z. Yuan, M. R. Hesamzadeh, Y. Cui, and L. B. Tjernberg, "Applying high performance computing to probabilistic convex optimal power flow," in *2016 International Conference on Probabilistic Methods Applied to Power Systems (PMAPS)*, Oct 2016, pp. 1–7.
[9] L. F. M. Machado, S. G. D. Santo, G. M. Junior, R. Itiki, and M. D. Manjrekar, "Multi-source distributed energy resources management system based on pattern search optimal solution using nonlinearized power flow constraints," *IEEE Access*, vol. 9, pp. 30 374–30 385, 2021.
[10] Q. Yu, Z. Jiang, Y. Liu, L. Li, and G. Long, "Optimization of an offshore oilfield multi-platform interconnected power system structure," *IEEE Access*, vol. 9, pp. 5128–5139, 2021.
[11] X. Bai, H. Wei, K. Fujisawa, and Y. Wang, "Semidefinite programming for optimal power flow problems," *International Journal of Electrical Power & Energy Systems*, vol. 30, no. 6, pp. 383–392, 2008.
[12] M. Farivar and S. H. Low, "Branch Flow Model: Relaxations and Convexification: Part I," *IEEE Transactions on Power Systems*, vol. 28, no. 3, pp. 2554–2564, Aug 2013.
[13] L. Gan, N. Li, U. Topcu, and S. H. Low, "Exact convex relaxation of optimal power flow in radial networks," *IEEE Transactions on Automatic Control*, vol. 60, no. 1, pp. 72–87, 2015.
[14] Z. Yuan and M. Paolone, "Properties of convex optimal power flow model based on power loss relaxation," *Electric Power Systems Research*, vol. 186, p. 106414, 2020.
[15] C. Chen, A. Atamtürk, and S. S. Oren, "Bound tightening for the alternating current optimal power flow problem," *IEEE Transactions on Power Systems*, vol. 31, no. 5, pp. 3729–3736, 2016.
[16] J. Bezanson, A. Edelman, S. Karpinski, and V. B. Shah, "Julia: A fresh approach to numerical computing," *SIAM Review*, vol. 59, no. 1, pp. 65–98, 2017.
[17] R. D. Zimmerman, C. E. Murillo-Sánchez, and R. J. Thomas, "Matpower: Steady-state operations, planning, and analysis tools for power systems research and education," *IEEE Transactions on Power Systems*, vol. 26, no. 1, pp. 12–19, 2011.
[18] J. H. Chow, *Time-Scale Modeling of Dynamic Networks with Applications to Power Systems*. Springer-Verlag, 1982.
[19] U. of Washington, "Power systems test case archive." [Online]. Available: http://www.ee.washington.edu/research/pstca/
[20] S. Fliscounakis, P. Panciatici, F. Capitanescu, and L. Wehenkel, "Contingency ranking with respect to overloads in very large power systems taking into account uncertainty, preventive, and corrective actions," *IEEE Transactions on Power Systems*, vol. 28, no. 4, pp. 4909–4917, Nov 2013.
[21] A. B. Birchfield, T. Xu, K. M. Gegner, K. S. Shetye, and T. J. Overbye, "Grid structural characteristics as validation criteria for synthetic networks," *IEEE Transactions on Power Systems*, vol. 32, no. 4, pp. 3258–3265, 2017.